\documentclass[12pt]{amsart}
\usepackage{graphicx}
\usepackage{amssymb, amsmath, amsthm, amsfonts, enumerate}
\usepackage{amsmath,setspace,scalefnt}
\usepackage{verbatim}
\usepackage{multirow}
\usepackage{mathtools}
\usepackage{comment}
\usepackage{float}
\usepackage{enumitem}
\restylefloat{figure}
\usepackage[margin=3.5cm]{geometry}

\usepackage{fancyhdr}
\usepackage{hyperref}

\newcommand{\Z}{{\mathbb Z}}
\newcommand{\N}{{\mathbb N}}

\newtheorem{theorem}{Theorem}
\newtheorem{lemma}[theorem]{Lemma}

\newtheorem{corollary}[theorem]{Corollary}

\newtheorem*{definition*}{Definition}

\newtheorem{remark}[theorem]{Remark}

\numberwithin{equation}{section}
\numberwithin{theorem}{section}

\newcommand{\zz}{\mathbb{Z}}

\newcommand{\set}[1]{\left\lbrace #1 \right\rbrace}

\subjclass[2010]{\ 05D10, 11B30, 11T06, 11A07}
\keywords{\ Ramsey theory, Arithmetic combinatorics}

\title[Towards characterizing the 2-Ramsey equations of the form $ax+by=p(z)$]
  {Towards characterizing the 2-Ramsey equations of the form $ax+by=p(z)$}

\author{Zsolt Baja}
\email{zsolt.baja@stud.ubbcluj.ro}
\address{Babeş-Bolyai University,
Faculty of Mathematics and Computer Sciences,
Kogălniceanu Street Nr. 1
400084 Cluj-Napoca, Romania}

\author{D\'aniel Dob\'ak}
\email{dd542@cam.ac.uk}
\address{University of Cambridge, Centre for Mathematical Sciences, Wilberforce Rd, Cambridge CB3 0WA, United Kingdom}

\author{Benedek Kov\'acs}
\email{benoke98@student.elte.hu}
\address{Institute of Mathematics, Faculty of Science, E\"otv\"os Lor\'and University, P\'azm\'any P\'eter s\'et\'any 1/C, H-1117 Budapest, Hungary.}

\author{P\'eter P\'al Pach}
\email{ppp@cs.bme.hu}
\address{Department of Computer Science and Information Theory, Budapest University of Technology and Economics, M\H{u}egyetem rkp. 3., H-1111 Budapest, Hungary; MTA-BME Lend\"ulet Arithmetic Combinatorics Research Group,
  ELKH, M\H{u}egyetem rkp. 3., H-1111 Budapest, Hungary.}

\author{Don\'at Pigler}
\email{pigler.donat@gmail.com}
\address{Institute of Mathematics, Faculty of Science, E\"otv\"os Lor\'and University, P\'azm\'any P\'eter s\'et\'any 1/C, H-1117 Budapest, Hungary.}

\begin{document}


\begin{abstract}

In this paper, we study a Ramsey-type problem for equations of the form $ax+by=p(z)$. We show that if certain technical assumptions hold, then any 2-colouring of the positive integers admits infinitely many monochromatic solutions to the equation $ax+by=p(z)$. 
This entails the $2$-Ramseyness of several notable cases such as the equation $ax+y=z^n$ for arbitrary $a\in\mathbb{Z}^{+}$ and $n\ge 2$, and also  
of $ax+by=a_Dz^D+\dots+a_1z\in\mathbb{Z}[z]$ such that $\textup{gcd}(a,b)=1$, $D\ge 2$, $a,b,a_D>0$ and $a_1\neq0$.
\end{abstract}

\maketitle

\section{Introduction}

The study of Ramsey theory searches for monochromatic patterns in finite colourings of $\zz^{+}$. It has a long history dating back to the famous theorem of Schur~\cite{Sch} in 1916, which states that the equation $x+y=z$ is \emph{Ramsey}, that is, any finite colouring of $\zz^{+}$ contains a monochromatic solution to $x+y=z$. Another classical example is van der Waerden’s theorem~\cite{vdW} stating that $\{x, x+y,\ldots, x+(\ell-1)y\}$ is Ramsey for any $\ell\in\zz^{+}$. Rado~\cite{Rad33} later in his seminal work resolved the Ramsey problem for all \emph{linear} equations, characterising all those that are Ramsey. Since then, many extensions have been studied, see e.g.~the polynomial extension of van der Waerden’s theorem by Bergelson and Leibman~\cite{BL96}.

In this paper, we study the polynomial extension of Schur's theorem, namely equations of the form $ax+by=p(z)$. We call an equation \emph{$k$-Ramsey}, $k\in\Z^+$, if any $k$-colouring of $\Z^+$ contains infinitely many monochromatic solutions to it.  Csikv\'ari, Gyarmati and S\'ark\"ozy~\cite{CGS} showed that $x+y=z^2$ is \emph{not} $16$-Ramsey, that is, they constructed a $16$-colouring of $\zz^{+}$ with no monochromatic solution for $x+y=z^2$ other than the trivial solution $x=y=z=2$. Green and Lindqvist~\cite{GL} completely resolved this case using Fourier-analytic arguments, giving the satisfying answer that any $2$-colouring of $\zz^{+}$ contains \emph{infinitely} many monochromatic solutions, and $3$ colours suffice to avoid non-trivial monochromatic solutions. In other words, $x+y=z^2$ is $2$-Ramsey, but not $3$-Ramsey. In fact, the $3$-colouring in~\cite{GL} can be easily adapted to show that $ax+by=p(z)$ is not $3$-Ramsey for any $a,b\in \mathbb{N}$ and $p(z)\in\Z[z]$ with $\deg(p)\ge 2$.
The fourth author of the present paper~\cite{Pach} gave a shorter combinatorial proof for the 2-Ramseyness of $x+y=z^2$. More recently, Liu, S\'andor and the fourth author of this paper~\cite{LPS} completely resolved the Ramsey problem for $$\{x,y,z:x+y=p(z)\}$$ 
for all polynomials over $\Z$, thus giving a polynomial extension of Schur's theorem. They proved that the equation $x+y=p(z)$ is 2-Ramsey if and only if $p(1)p(2)$ is even. If $p(1)p(2)$ is odd, then $p$ attains only odd values, thus colouring the integers according their parity avoids monochromatic solutions to $x+y=p(z)$. According to their result this divisibility barrier is the only obstruction to the 2-Ramseyness for $x+y=p(z)$. They also ask what happens for other linear forms in place of $x+y$.

We study the Ramsey problem for equations of the form $ax+by=p(z)$, where $a$ and $b$ are positive integers. As mentioned earlier, using the greedy colouring it can be seen that the equation is {\it not} 3-Ramsey. In the case $a=b=1$ the equation turned out to be 2-Ramsey unless the odd-even colouring avoids monochromatic solutions. On the other hand, if coefficients larger than 1 are allowed, then 2-Ramseyness can fail to hold in several cases. Before describing our results in the positive direction let us show a few examples when the equation is not 2-Ramsey. Throughout the paper our colourings (usually denoted by $\phi$) will use the two colours $+1$ and $-1$, which we will often shorten to $+$ and $-$.

\medskip
\noindent
{\emph {Example} 1.} If for every $z$ we have $\gcd (a,b)\nmid p(z)$, then $ax+by=p(z)$ is trivially not 2-Ramsey, in fact the equation does not have an integer solution at all.

\medskip
\noindent
{\emph {Example} 2.} Let $p$ be a prime, then the equation $ax+by=z^{p-1}-1$ is not 2-Ramsey if $p\mid a,\ p\nmid b$. Indeed, let $\phi(t)=+$, if $p\mid t$ and $\phi(t)=-$, otherwise. If $p\mid x,y,z$, then $ax+by\equiv 0\not\equiv -1\equiv z^{p-1}-1 \pmod{p}$. If $p\nmid x,y,z$, then  $ax+by\not\equiv 0\equiv z^{p-1}-1 \pmod{p}$.

More generally, for any $n\geq 1$ the equation $ax+by=z^{p^n-p^{n-1}}-1$ is not 2-Ramsey if $p^n\mid a,\ p^n\nmid b$. 

\medskip
\noindent
{\emph {Example} 3.} The equation $x+4y=(z+2)(z+1)=z^2+3z+2$ is not 2-Ramsey. Indeed, for the colouring $\phi(t)=+$ for $t\equiv 0,1\pmod{4}$, $\phi(t)=-$ for $t\equiv 2,3\pmod{4}$ no monochromatic solution exists.

\medskip

Our main results are the following:

\begin{theorem}\label{cor-const0}
Let $a$ and $b$ be positive integers such that $\gcd(a,b)=1$. If $p(z)=a_Dz^D+\dots+a_1z\in \mathbb{Z}[z]$ satisfies $D\ge 2$, $a_D>0$ and $a_1\neq 0$, then
the equation $ax+by=p(z)$ is 2-Ramsey.
\end{theorem}

\begin{theorem}\label{cor-b=1}
Let $a$ and $n$ be positive integers with $n\ge 2$. Then 
the equation $ax+y=z^n$ is 2-Ramsey.
\end{theorem}

The proofs of Theorem~\ref{cor-const0} and Theorem~\ref{cor-b=1} go along similar lines. In fact both statements follow from the following more general theorem which shows that when certain technical conditions hold, the equation $ax+by=p(z)$ is 2-Ramsey:

\begin{theorem}
\label{mainresult}
Suppose that $p\in \mathbb{Z}[z]$ is a polynomial of degree at least 2 with a positive leading coefficient, and let $a$ and $b$ be positive integers with $\gcd(a,b)=1$. Also suppose that there exist positive integers $d$ and $u$, and integers $t$ and $v$ (which will represent residue classes mod $d$ and mod $u$ respectively) such that the following six conditions hold:

\begin{enumerate}
  \item $(a+b)t\equiv p(t) \pmod{d}$,
  \item if  $k\equiv t\pmod{d}$, then $a\mid p(k+d)-p(k)$,
  \item if  $k\equiv t\pmod{d}$, then letting $m=\frac{p(k+d)-p(k)}{a}$ and $m'=\frac{m}{u}$, we have $u\mid m$ and $2\nmid m'$ and $\gcd(a,m')=\gcd(b,m')=\gcd(u,m')=1$,
  \item $p(v)\equiv (a+b)v~\pmod{u}$,
  \item if  $k\equiv t\pmod{d}$, then $p(k)\equiv (a+b)v\pmod{u}$,
  \item if  $j\equiv v\pmod{u}$ and $k\equiv t\pmod{d}$, then $b\mid p(k)-aj$.
\end{enumerate}

Then the equation $ax+by=p(z)$ is 2-Ramsey.
\end{theorem}

\begin{remark}
Note that if the hypotheses of Theorem~\ref{mainresult} are satisfied for $p(z)$ with some parameters $d,u,t,v$, then they are also satisfied for $p(z)+\lambda$ (with the same parameters) for any $\lambda$ which is simultaneously divisible by $d,u$ and $b$. 
\end{remark}


\noindent
Some further applications of Theorem~\ref{mainresult} are as follows.

\begin{corollary}\label{cor-czp}
Let $a,b,c$ be positive integers with $\gcd(a,b)=1$ and let $p$ be a prime such that $p\nmid a$. Write $b=p^{pn-\varepsilon}B$, where $\nu_p(b)=pn-\varepsilon$, $n\in\zz$ and $\varepsilon\in\{0,1,\ldots,p-1\}$. If $\gcd(a,c)=1$, $\gcd(B,c)=1$ and there exists $t\in \zz_{aBp^n}$ such that
\[
    \left\lbrace{
    \begin{array}{ll}
        t \equiv0 &\pmod{p^n}\\
    t^{p-1} \equiv ac^{-1} &\pmod{B}\\
    t^{p-1} \equiv bc^{-1} &\pmod{a}
    \end{array}
    }\right.\]
then the equation $ax+by=cz^p$ is 2-Ramsey.
\end{corollary}

Note that in the special case $p=2$, Corollary~\ref{cor-czp} yields the following:

\begin{corollary}\label{rem-cz2}
Let $a,b,c$ be positive integers such that $\gcd(a,b)=1$, $\gcd(a,c)=2^\ell$ and $\gcd(b,c)=2^{\ell'}$ for some $\ell,\ell'\in\zz^{\ge 0}$.
Then the equation $ax+by=cz^2$ is 2-Ramsey.
\end{corollary}

\begin{corollary}\label{cor-dcz2}
Let $a$, $b$, $c$ be positive integers such that $\gcd\left(c,\frac{a}{\gcd(a,b)}\right)=2^\ell$ and $\gcd\left(c,\frac{b}{\gcd(a,b)}\right)=2^{\ell'}$ for some $\ell,\ell'\in\zz^{\ge 0}$. Then the equation $ax+by=cz^2$ is 2-Ramsey.
\end{corollary}

Note that a special case of Corollary~\ref{cor-dcz2} is that the equation $ax+by=z^2$ is 2-Ramsey for every $a,b\in\mathbb{Z}^+$.

\subsection{Other related work}

It is worth noting that the Ramsey problem for $x^{\alpha}+y^{\beta}=z^{\gamma}$ in $\Z/p\Z$ has been studied by Lindqvist~\cite{Lin}. If one puts no restriction on $z$ and looks for monochromatic pairs $\{x,y\}$ with $x+y$ being a perfect square, then Khalafallah and Szemer\'edi~\cite{KS} showed that this is Ramsey in $\mathbb{Z}^+$. Yet another similar looking pattern that behaves very differently is to consider $x-y$ instead. Bergelson~\cite{Ber}, improving upon results of Furstenberg~\cite{Fur77} and S\'ark\"ozy~\cite{Sar78}, proved that $\{x,y,z:x-y=z^2\}$ is Ramsey.

Ramsey theory has witnessed exciting development recently. We refer the readers to the papers of Green and Sanders~\cite{GS16} and of Moreira~\cite{Mor17} for the problem involving the sum and the product of $x$ and $y$, and to the papers of Di Nasso and Luperi Baglini~\cite{NB18}, and of Chow, Lindqvist and Prendiville~\cite{CLP18+} for  generalisations of Rado's criterion to non-linear polynomials.

\subsection*{Notations.} Whenever $m$ is a positive integer, we shall write $\zz_m:=\zz/m\zz$. If $p$ is a prime and $a$ is a nonzero integer, then $\nu_p(a)$ denotes the largest integer such that $p^{\nu_p(a)}\mid a$. 

\subsection*{Organization. } The rest of the paper is organised as follows. In Section~\ref{sec-2} we prove Theorem~\ref{mainresult}, in Section~\ref{sec-3} we present the proofs of Theorem~\ref{cor-const0}, Theorem~\ref{cor-b=1}, Corollary~\ref{cor-czp}, Corollary~\ref{rem-cz2} and Corollary~\ref{cor-dcz2}.

\section{Proof of Theorem~\ref{mainresult}}
\label{sec-2}


\noindent
We first prove some technical lemmas necessary for the proof. 

\begin{lemma}\label{bigsolinresidueclass}
Let $p\in \zz[z]$ be a polynomial of degree at least 2 with a positive leading coefficient, and let $a,b\in \zz^{+}$ with $\gcd(a,b)=1$. If $d\ge 1$ is an integer and $d\zz+t$ is a residue class such that $(a+b)t\equiv p(t)\pmod{d}$, then for any $K\in \zz^{+}$, there exists a solution $(x,y,z)$ of $ax+by=p(z)$ such that $x,y,z\ge K$ and $x,y,z\in d\zz+t$.
\end{lemma}
\begin{proof}
Since $p$ is of degree at least 2 and its leading coefficient is positive, we can choose a value $z_1\ge K$ with $z_1\equiv t\pmod{d}$ such that $p(z_1)\ge (a+b)z_1+4a^2b^2d$. Now 
$$p(z_1)\equiv p(t)\equiv (a+b)t\equiv (a+b)z_1\pmod{d},$$
say, $p(z_1)=(a+b)z_1+\ell d$, where we know that $\ell\ge 4a^2b^2$. Let us write $\ell=\ell_1\cdot 2ab+\ell_2$, where $\ell_1, \ell_2$ are integers with $0\le \ell_2<2ab$. The lower bound on $\ell$ gives $\ell_1\ge 2ab$.

Since $\gcd(a,b)=1$, there exist integers $r$ and $s$ such that $ra+sb=1$. It is easy to see that we can take $r$ and $s$ with $|r|\le b$ and $|s|\le a$. Let us consider
$$x=z_1+(\ell_1b+\ell_2r)d,$$
$$y=z_1+(\ell_1a+\ell_2s)d,$$
$$z=z_1.$$

Now, as $|\ell_2|< 2ab$ and $|r|\le b$, using the fact that $\ell_1\ge 2ab$ we have $\ell_1b+\ell_2r\ge 0$. Similarly $\ell_1a+\ell_2s\ge 0$ as well, so $x,y,z\ge z_1\ge K$ with $x,y,z\equiv t\pmod{d}$. Furthermore, 
\begin{multline*}
ax+by=(a+b)z_1+(2ab\ell_1+(ra+sb)\ell_2)d=(a+b)z_1+(2ab\ell_1+\ell_2)d=\\
=(a+b)z_1+\ell d=p(z_1)=p(z),
\end{multline*}

so we have found a solution satisfying the required properties.
\end{proof}

\begin{lemma}\label{polyvalue_inlineargap}
Let $p\in \zz[z]$ be a polynomial of degree at least 1 with a positive leading coefficient. Let $d\in \zz^{+}$, take any residue class $d\zz+t$ and any real number $\delta>0$. Then for sufficiently large real values of $x$ there exists $z\in d\zz+t$ such that $x\le p(z)\le (1+\delta) x$.
\end{lemma}
\begin{proof}
Let $p(z)=c_nz^n+c_{n-1}z^{n-1}+\dots+c_1z+c_0$. For sufficiently large $z'$ we have $p(d(z'+1)+t)-p(dz'+t)\le \delta p(dz'+t)$, since the degree in $z'$ is $n-1$ on the left hand side and $n$ on the right hand side, with positive leading coefficients ($nc_nd^n$ and $\delta c_nd^n$ respectively). Choose $z'_0$ to be an integer so that this holds for $z'\ge z'_0$, while also choosing $z'_0$ to be large enough so that $p$ is increasing on the interval $[dz'_0+t, \infty)$. Then we claim that the statement of the lemma holds for all $x>p(dz'_0+t)$. Fixing such a value $x$, take $z'$ to be the largest integer such that $p(dz'+t)<x$. Then we have $z'\ge z'_0$, and taking $z=d(z'+1)+t$, we have $x\le p(z)\le (1+\delta)p(dz'+t)<(1+\delta)x$. 
\end{proof}

\begin{lemma}\label{resclass_bijection_lemma}
 Let $u,m',a_1,a_2\in \zz^{+}$ and $C,v\in \zz$ such that $\gcd(u,m')=\gcd(a_1,m')=\gcd(a_2,m')=1$ and $(a_1+a_2)v\equiv C\pmod{u}$. Then denoting $m=um'$, the following statement holds: for every $\gamma_1\in u\zz_m+v$, there exists a unique $\gamma_2\in u\zz_m+v$ such that $a_1\gamma_1+a_2\gamma_2\equiv C\pmod{m}$, and the mapping $u\zz_m+v\to u\zz_m+v$, $\gamma_1\mapsto \gamma_2$ is a bijection.
\end{lemma}
\begin{proof}
Since for all $\gamma_1,\gamma_2\in u\zz_m+v$ we have $a_1\gamma_1+a_2\gamma_2\equiv (a_1+a_2)v\equiv C\pmod{u}$, and $m=um'$ with $\gcd(u,m')=1$, the condition $a_1\gamma_1+a_2\gamma_2\equiv C\pmod{m}$ is equivalent to $a_1\gamma_1+a_2\gamma_2\equiv C\pmod{m'}$.

Observe that the representatives $v, u+v, 2u+v, \dots, (m'-1)u+v$ of the elements of $u\zz_m+v$ form a complete residue system mod $m'$ since $\gcd(u,m')=1$. Since $\gcd(a_2,m')=1$, these representatives still form a complete residue system mod $m'$ once all of them are multiplied by $a_2$, so for a given $\gamma_1\in u\zz_m+v$, there will be a unique $\gamma_2$ such that $a_2\gamma_2\equiv C-a_1\gamma_1\pmod{m'}$. The roles of the indices $1$ and $2$ can be reversed in this argument, giving that our mapping is bijective.
\end{proof}

\begin{proof}[Proof of Theorem~\ref{mainresult}]

First we quickly check the statement in the case $a=b=1$ by using the necessary and sufficient condition given in \cite[Corollary 1.4]{LPS}. According to this result, for a polynomial $p$ of degree at least 1 and with positive leading coefficient, $x+y=p(z)$ is 2-Ramsey if and only if $2\mid p(1)p(2)$. So it suffices to show that if $p(z)\equiv 1\pmod{2}$ for all integers $z$, then the conditions (1)-(6) cannot all hold. Indeed, if our conditions all hold for $d,u,t,v$, then (4) implies that $u\mid p(v)-2v$, and as $p(v)$ is odd, $u$ must also be odd. For $k=t$, we have $m=p(k+d)-p(k)$ which is necessarily even, and as $u$ is odd, $m'=\frac{m}{u}$ (which is an integer by (3)) is also even. However, this contradicts that by (3), $m'$ should be odd. This completes the proof of the case $a=b=1$. 

From now on, we may assume $a$ and $b$ are distinct. Let us denote $U=\max(a,b)$ and $L=\min(a,b)$. Let $p(z)=c_nz^n+c_{n-1}z^{n-1}+\dots+c_1z+c_0$, where $n\ge 2$ and $c_n>0$. Let us take $d,u\in \zz^{+}$ and residue classes $d\zz+t$ and $u\zz+v$ that satisfy the conditions (1)-(6). Let  $\phi: \zz^{+}\to \{+,-\}$ be an arbitrary 2-colouring, our aim is to show that there exist infinitely many monochromatic solutions to $ax+by=p(z)$. For the sake of contradiction, assume that there are only finitely many.

In our proof, we will use real constants $\varepsilon_0>0$ and $\varepsilon_1>0$ chosen to be sufficiently small. We let $\varepsilon_0<\frac{1}{U+1}$, whereas the value of $\varepsilon_1$ will be specified later depending on $a$, $b$, $p$ and $d$ (according to condition $(*)$ in the proof of Claim 4).

For a positive integer $k$ with $k\equiv t\pmod{d}$, we say that a \textit{$d$-switch} occurs at $k$ if $\phi(k)\ne \phi(k+d)$. Let us call two $d$-switches $k_1<k_2$ \textit{neighbouring} if there is no $k_1<k_3<k_2$ with a $d$-switch at $k_3$.

By Lemma \ref{bigsolinresidueclass}, we can see that there must be infinitely many $d$-switches: if the set $(d\zz+t)\cap [K,\infty)$ was monochromatic for some $K\in \zz^{+}$, then any solution $(x,y,z)$ in this set (which exists by the Lemma) would be monochromatic. We note that the condition $(a+b)t\equiv p(t)\pmod{d}$ of the Lemma is granted by (1). By increasing $K$ toward infinity, we would get infinitely many monochromatic solutions for a contradiction.

\bigskip

\textbf{Claim 1.} There are infinitely many pairs $k_1<k_2$ of neighbouring $d$-switches such that $k_2<Uk_1$.

\smallskip

\textbf{Proof.} Assume for a contradiction that there are only finitely many such pairs of $d$-switches, with $K_0$ being an upper bound for all of them. Let us choose an integer $K$ large enough, so that $K>K_0$, and $K$ also fulfills some other conditions implied in the later description.

Take neighbouring $d$-switches $k_1<k_2<k_3$ such that $k_1>K$. Then as $k_1,k_2,k_3>K_0$, we have $k_2\ge Uk_1$ and $k_3\ge Uk_2$. Without the loss of generality we may assume that $\phi(x)=-$ for all numbers $x\in (d\zz+t)\cap (k_1,k_2]$ and $\phi(x)=+$ for all numbers $x\in (d\zz+t)\cap (k_2,k_3]$. For technical reasons, fix a small $\varepsilon>0$ such that 
$L+\frac{L}{U}+\varepsilon<U-\varepsilon$ (this is possible, since $U\ge L+1$ implies that $L+\frac{L}{U}<U$).

Now, using Lemma \ref{polyvalue_inlineargap} and assuming that $K$ (and hence $k_2$) is sufficiently large, let us pick $z\in d\zz+t$ such that $$\left(L+\frac{L}{U}+\varepsilon\right)k_2\le p(z)\le (U-\varepsilon)k_2.$$

Let us distinguish two cases depending on the colour of $z$.

\medskip

\textit{Case 1: $\phi(z)=-$.}

In this case, our aim is to find positive integers $x,y\in d\zz+t$ such that $Lx+Uy=p(z)$ and $k_1<x,y\le k_2$. Then we will have $\phi(x)=\phi(y)=-$, which will give a monochromatic solution with $\max(x,y,z)\ge K$.

First, we prove that there exist $x_1,y_1\in d\zz+t$ with $Lx_1+Uy_1=p(z)$. Let us take $r,s\in \zz$ such that $rL+sU=1$. Then by (1), 
$$p(z)\equiv (L+U)t\pmod{d},$$
so $p(z)-(L+U)t=d\ell$ for some $\ell\in \zz$. Now let $x_1=t+r\ell d$ and $y_1=t+s\ell d$, so then 
$$Lx_1+Uy_1=L(t+r\ell d)+U(t+s\ell d)=(L+U)t+\ell d(rL+sU)=(L+U)t+\ell d=p(z).$$

Now, observe that for any $w\in \zz$, the numbers $x=x_1+wdU$ and $y=y_1-wdL$ also satisfy $Lx+Uy=p(z)$ and $x,y\in d\zz+t$. Every time we add 1 to $w$, the quantity $x-y$ increases by $d(L+U)$. So there exists $w\in \zz$ such that $|x-y|\le d\cdot \frac{L+U}{2}$, and so by considering the point where $x$ and $y$ would meet if they moved continuously, $$\frac{p(z)}{L+U}-d\cdot \frac{L+U}2\le x,y\le \frac{p(z)}{L+U}+d\cdot \frac{L+U}2.$$
Then for $K$ (and hence $k_2$) being sufficiently large,
$$x,y\ge \frac{\left(L+\frac{L}{U}+\varepsilon\right)k_2}{L+U}-d\cdot \frac{L+U}{2}=\frac{LU+L}{U(L+U)}k_2+\frac{\varepsilon}{L+U}k_2-d\cdot \frac{L+U}{2}>\frac{1}{U}k_2\ge k_1$$
and
$$x,y\le \frac{U-\varepsilon}{L+U}k_2+d\cdot \frac{L+U}2<k_2$$
where we used that $\frac{1}{U}\le \frac{LU+L}{U(L+U)}$. So we have found $x,y$ satisfying our aim.

\medskip

\textit{Case 2: $\phi(z)=+$.}

As in the previous case, there exist $x_1,y_1\in d\zz+t$ such that $Lx_1+Uy_1=p(z)$. Again, taking $x=x_1+wdU$ and $y=y_1-wdL$ for some $w\in \zz$, now our goal is to obtain $x,y>0$ with $\phi(x)=\phi(y)=+$.

Since we have assumed that there are only finitely many ``close'' pairs of $d$-switches, it follows that there must exist $L$ consecutive elements of $d\zz+t$ such that all of them have the colour $+$ (otherwise there would be infinitely many pairs of neighbouring switches with a difference of at most $(L-1)d$, which would eventually be ``close'' pairs). It is possible to set $w$ such that $y=y_1-wdL$ is equal to one of these. Then 
$$x=\frac{p(z)}{L}-\frac{U}{L}y\ge \left(1+\frac{1}{U}+\frac{\varepsilon}{L}\right)k_2-\frac{U}{L}y>k_2$$
for $K$ (and hence $k_2$) being sufficiently large. Also, clearly $x\le p(z)<Uk_2\le k_3$. Thus $\phi(x)=+$ also holds, since $x\in (k_2,k_3]$. Therefore, we found again a monochromatic solution with $\max(x,y,z)\ge K$ (since $x>k_2>K$).

~

In both cases we found a monochromatic solution with $\max(x,y,z)\ge K$. Since $K$ was arbitrarily large, this means that $ax+by=p(z)$ has infinitely many monochromatic solutions; finishing off the proof of Claim 1. $\qed$

~

In the following arguments, let $k_1$ and $k_2$ be neighbouring $d$-switches with $k_1<k_2<Uk_1$. By Claim 1, $k_1$ can be arbitrarily large; let us take it to be larger than $\frac{1}{\varepsilon_0}M$ where $M$ is the maximal value of $\max(x,y,z)$ over all monochromatic solutions to $ax+by=p(z)$. (We have set $\varepsilon_0<\frac{1}{U+1}<1$ so this is also greater than $M$.) Later we will also specify some additional conditions on how large $k_1$ will be taken. Now, let $k$ be either $k_1$ or $k_2$. Without loss of generality we may assume that $\phi(k)=+$ and $\phi(k+d)=-$.

Let us take any integer $0<j<\frac{p(k)}{a}$ such that $j\equiv v\pmod{u}$. Then by (6), we have $b\mid p(k)-aj$. Consider the following two equations where $m=\frac{p(k+d)-p(k)}{a}$, as defined in condition (3) (an integer by (2)):
$$aj+b\cdot \frac{p(k)-aj}{b}=p(k),$$
$$a(j+m)+b\cdot \frac{p(k)-aj}{b}=p(k+d).$$
As all monochromatic solutions to $ax+by=p(z)$ have $\max(x,y,z)<k$, none of the triples $\left(j,\frac{p(k)-aj}{b},k\right)$ and $\left(j+m, \frac{p(k)-aj}{b}, k+d\right)$ can be monochromatic. This yields that $\phi(j)\le \phi(j+m)$, since otherwise $\phi(j)=+$ and $\phi(j+m)=-$, and no matter the colour of $\frac{p(k)-aj}{b}$, one of the two triples mentioned will be monochromatic.

So for each congruence class $\gamma\in \zz_m$ such that $\gamma\equiv v\pmod{u}$ (note that this makes sense as $u\mid m$), $\phi$ is monotonic on the members of $(m\zz+\gamma)\cap\left[1, \frac{p(k)}{a}\right)$. For every such congruence class $\gamma$, let $\beta(\gamma)$ denote the smallest member of $(m\zz+\gamma)\cap\left[1, \frac{p(k)}{a}\right)$ that is coloured $+$. If there is no such member, let $\beta(\gamma)=\frac{p(k)}{a}$.

Now, let us define the following set:
$$A=\left\{\gamma\in u\zz_m+v:\  \beta(\gamma)\ge \frac{p(k)}{a+b+\varepsilon_1}\right\}.$$

Now we use Lemma \ref{resclass_bijection_lemma} to assign to every element $\gamma\in u\zz_m+v$ an element $\gamma'\in u\zz_m+v$ such that $a\gamma'\equiv p(k)-b\gamma\pmod{m}$, giving a bijection $u\zz_m+v\to u\zz_m+v, \gamma\mapsto \gamma'$. The Lemma is used with $u$, $v$ and $m'=\frac{m}{u}$ (as in condition (3)), and $a_1=b$, $a_2=a$ and $C=p(k)$. The required conditions are $\gcd(u,m')=\gcd(a,m')=\gcd(b,m')=1$ and $(a+b)v\equiv p(k)\pmod{u}$, guaranteed by (3) and (5) respectively.

\bigskip

\textbf{Claim 2.} For every $\gamma\in u\zz_m+v$, either $\gamma$ or $\gamma'$ belongs to $A$.

\smallskip

\textbf{Proof.} Suppose neither of them does. In this case we will find a monochromatic solution to $ax+by=p(z)$ with $z=k$, which gives a contradiction (as we know that all such solutions have $\max(x,y,z)<k_1$). Take $r,s\in \zz$ such that $ra+sb=1$, and pick some integers $x_0\equiv \gamma'$ and $y_0\equiv \gamma\pmod{m}$. As we know that $a\gamma'+b\gamma\equiv p(k)\pmod{m}$, we can write $ax_0+by_0=p(k)+\ell m$ for some $\ell\in \zz$. Taking $x_1=x_0-r\ell m$ and $y_1=y_0-s\ell m$, we have $ax_1+by_1=p(k)$. For any $w\in \zz$, if we let $x=x_1+wbm$ and $y=y_1-wam$, then we have $ax+by=p(k)$ with $x\equiv \gamma' \pmod{m}$ and $y\equiv \gamma\pmod{m}$. Similarly to the argument seen in the proof of Claim 1, there exists some $w$ such that 
$$\frac{p(k)}{a+b}-\frac{a+b}{2}\cdot m\le x,y\le \frac{p(k)}{a+b}+\frac{a+b}{2}\cdot m.$$
As $\varepsilon_1$ will be chosen (in the proof of Claim 4) independently from $k_1$, and $m$ is a polynomial in $k$ having smaller degree than $p$, we can take $k_1$ (and hence $k$) large enough so that $$\frac{p(k)}{a+b+\varepsilon_1}\le\frac{p(k)}{a+b}-\frac{a+b}{2}\cdot m\le x,y\le \frac{p(k)}{a+b}+\frac{a+b}{2}\cdot m<\frac{p(k)}{a}.$$

Now as $\gamma,\gamma'\not\in A$, all members of $m\zz+\gamma$ and $m\zz+\gamma'$ lying in $\left[\frac{p(k)}{a+b+\varepsilon_1}, \frac{p(k)}{a}\right)$ are coloured $+$, since by the definition of $\beta(\gamma)$, every element of $(m\mathbb{Z}+ \gamma)\cap \left[\beta(\gamma),\frac{p(k)}{a}\right)$ is coloured $+$. Hence, $\phi(x)=\phi(y)=+$, and we also know that $\phi(k)=+$, so with $z=k$ we have a monochromatic solution $(x,y,z)$, contradiction. $\qed$

\bigskip

\textbf{Claim 3.} We have $aA+bA=u\zz_m+p(v)$.

\smallskip

\textbf{Proof.} Consider the permutation $i: \gamma\mapsto \gamma'$ on $u\zz_m+v$, and write $i$ as a product of disjoint cycles. By Claim 2, in each cycle, if we consider any two neighbouring elements, then at least one must belong to $A$, so at least half of the members of the cycle must be in $A$. Combined with the fact that the base set has odd cardinality ($m'$ is odd by (3)), we must have $|A|>\frac12|u\zz_m+v|$.

If we take elements $\gamma_1,\gamma_2\in A$, then $a\gamma_1+b\gamma_2\equiv (a+b)v\pmod{u}$, so by (4), $a\gamma_1+b\gamma_2\equiv p(v)\pmod{u}$. So certainly, $aA+bA\subseteq u\zz_m+p(v)$. Now, fix a residue $\delta\in u\zz_m+p(v)$. Again using Lemma \ref{resclass_bijection_lemma} but this time with $C=\delta$, we can show that there is a bijection $j: u\zz_m+v\to u\zz_m+v$ such that for every $\gamma$, we have $aj(\gamma)+b\gamma\equiv \delta\pmod{m}$. This is because we know that $\delta\equiv p(v)\equiv (a+b)v\pmod{u}$ by (4). This permutation can again be expressed as a product of disjoint cycles, and since $|A|>\frac12|u\zz_m+v|$, at least one cycle has more than half of its elements in $A$. This means that this cycle must have two adjacent elements in $A$, giving an element $\gamma\in A$ with $j(\gamma)\in A$. Then $\delta \equiv aj(\gamma)+b\gamma\pmod{m}$, so $\delta\in aA+bA$. This concludes the proof of Claim 3. $\qed$

\bigskip

\textbf{Claim 4.} All members of the set $[\varepsilon_0k,(1-\varepsilon_0)k]\cap (u\zz+v)$ have the colour $+$.

\smallskip

\textbf{Proof.} We have fixed $0<\varepsilon_0<\frac{1}{U+1}$ which means that $\frac{1}{U}(1-\varepsilon_0)>\varepsilon_0$. (Note that since $U\ge 2$, this also implies that $\varepsilon_0<\frac12$, so the given interval is meaningful.)

Take an integer $i\in [\varepsilon_0k,(1-\varepsilon_0)k]$ with $i\equiv v\pmod{u}$. Then $p(i)\equiv p(v)\pmod{u}$. In Claim 3 we have seen that no matter what the residue class of $p(i)$ is modulo $m$ (denote this class by $\overline{p(i)}$), we can write it as $\overline{p(i)}=a\overline{j}_1+b\overline{j}_2$ where $\overline{j}_1,\overline{j}_2\in A$. Fix any representatives $j_1$ and $j_2$ of the modulo $m$ classes $\overline{j}_1$ and $\overline{j}_2$ respectively. Then $aj_1+bj_2=p(i)+\ell m$ for some $\ell\in \zz$. Taking $ra+sb=1$, and replacing $j_1$ and $j_2$ by $j_1-r\ell m$ and $j_2-s\ell m$ respectively, we can assume that $\ell=0$, and hence $aj_1+bj_2=p(i)$.

We will show that there exist representatives $j'_1$ and $j'_2$ of the same classes such that $aj'_1+bj'_2=p(i)$ and $0<j'_1,j'_2<\frac{p(k)}{a+b+\varepsilon_1}$. If we manage to show this, then by the definition of $A$ we would have $\phi(j'_1)=\phi(j'_2)=-$, and based on the  hypothesis that $\phi(i)=-$, we would have a monochromatic solution to $ax+by=p(z)$ with $\max(x,y,z)\ge i\ge \varepsilon_0k>M$, which  contradicts the definition of $M$.

For $w\in \zz$, take $j'_1=j_1+wbm$ and $j_2'=j_2-wam$. Then $aj_1'+bj_2'=p(i)$, and by our usual argument, we can choose $w$ such that 
$$\frac{p(i)}{a+b}-\frac{a+b}{2}m\le j_1', j_2'\le \frac{p(i)}{a+b}+\frac{a+b}{2}m.$$

As $p(z)$ and $p(z+d)-p(z)$ are polynomials of degree $n$ and $n-1$ respectively with positive leading coefficients, there exist positive reals $K_1,K_2,C,D,E,F$ (depending on $p$, $d$ and $\varepsilon_0$) such that $Ez^n\le p(z)\le Fz^n$ for all reals $z\ge K_1$, and $Cz^{n-1}\le p(z+d)-p(z)\le Dz^{n-1}$ for all reals $z\ge K_2$. Here $\frac{F}{E}$ can be any real number greater than 1, so we can demand that $F<E\cdot \frac{1}{(1-\varepsilon_0)^n}$. Also there is a $K_3>0$ such that $p$ is increasing on $[K_3, \infty)$. 

Now, that we have $F(1-\varepsilon_0)^n<E$, we can pick $\varepsilon_1>0$ small enough such that
$$\frac{F(1-\varepsilon_0)^n}{a+b}<\frac{E}{a+b+\varepsilon_1}. \eqno(*)$$

For sufficiently large $k$, we have $\varepsilon_0k\ge \max(K_3,K_1)$, meaning that 
$$E(\varepsilon_0k)^n\le p(\varepsilon_0k)\le p(i) \le p((1-\varepsilon_0)k)\le F((1-\varepsilon_0)k)^n. \eqno(**)$$

Also for $k\ge K_2$ we have
$$\frac{Ck^{n-1}}{a}\le m=\frac{p(k+d)-p(k)}{a}\le \frac{Dk^{n-1}}{a}. \eqno({*}{*}{*})$$

Using $(*)$, $(**)$ and $({*}{*}{*})$, for sufficiently large $k$ we have

$$j_1',j_2'\le \frac{p(i)}{a+b}+\frac{a+b}{2}m\le \frac{F(1-\varepsilon_0)^n}{a+b}k^n+\frac{a+b}{2}\cdot \frac{D}{a}k^{n-1}<\frac{E}{a+b+\varepsilon_1}k^n$$

where for $k\ge K_1$ we have $Ek^n\le p(k)$, so $j_1',j_2'<\frac{p(k)}{a+b+\varepsilon_1}$.

Also, using $(**)$ and $({*}{*}{*})$, for sufficiently large $k$ we have

$$j_1',j_2'\ge \frac{p(i)}{a+b}-\frac{a+b}{2}m\ge \frac{E\varepsilon_0^n}{a+b}k^n-\frac{a+b}2\cdot \frac{D}{a}k^{n-1}>0$$

So we have shown that $\phi(i)=+$, concluding the proof of Claim 4. $\qed$

~

To finish off the proof of Theorem~\ref{mainresult}, observe that the argument so far can be repeated with both $k=k_1$ and $k=k_2$. Without loss of generality we can assume that $\phi(k_1)=+$ and $\phi(k_2)=-$, then (as we assumed $\phi(k)=+$) by Claim 4 we have that all members of $[\varepsilon_0k_1,(1-\varepsilon_0)k_1]\cap (u\zz+v)$ have the colour $+$, and all members of $[\varepsilon_0k_2, (1-\varepsilon_0)k_2]\cap (u\zz+v)$ have the colour $-$. If the two sets have a common member, then we get a contradiction. For this, a sufficient condition is that $(1-\varepsilon_0)k_1\ge \varepsilon_0k_2+u$. Since we have $\varepsilon_0<\frac{1}{U}(1-\varepsilon_0)$, for $k_1$ (and hence $k_2$) large enough, we have $\varepsilon_0k_2+u<\frac{1}{U}(1-\varepsilon_0)k_2<(1-\varepsilon_0)k_1$. (This is the crucial part where we used the fact that $k_2<Uk_1$.)

This concludes the proof of Theorem \ref{mainresult}. \end{proof}

\section{Proofs}
\label{sec-3}

\noindent
In this section we prove Theorem~\ref{cor-const0}, Theorem~\ref{cor-b=1} and the corollaries by using Theorem~\ref{mainresult}.

\begin{proof}[Proof of Theorem~\ref{cor-const0}]
Let us write $a=2^{\alpha_1}\prod\limits_{i=2}^r p_i^{\alpha_i}$ and  $b=\prod\limits_{j=1}^s q_j^{\gamma_j}$, where the $p_i$ and $q_j$ are pairwise distinct primes, $p_1=2$, $\alpha_1\ge 0$, and all other exponents are positive. (If $a$ is odd and $b$ is even, then we swap them.)

Let 
$$d=\prod_{i=1}^r p_i^{\max\{\alpha_i,\ \nu_{p_i}(a_1)+1\}}\cdot \prod_{j=1}^s q_j^{\max\{\gamma_j,\ \nu_{q_j}(a_1)+1\}}$$
and $t=0$ (note that as $a_1\neq 0$, $\nu_q(a_1)$ is well-defined for any prime $q$). Let 
$$u= \prod_{i=1}^r p_i^{\max\{\nu_{p_i}(a_1),\ 2\nu_{p_i}(a_1)+1-\alpha_i\}}\cdot \prod_{j=1}^s q_j^{\max\{\gamma_j+\nu_{q_j}(a_1),\ 2\nu_{q_j}(a_1)+1\}}$$
and $v=0$. These choices satisfy the conditions of Theorem~\ref{mainresult}, since, if we assume $d\mid k$ and $u\mid j$:

\begin{enumerate}
    \item $(a+b)\cdot 0 = 0 \equiv p(0) \pmod{d}$
    \item Since $a\mid d$, we have $a\mid p(k+d)-p(k)$.
    \item We have 
    $$m=\frac{p(k+d)-p(k)}{a}=\frac{\sum\limits_{i=1}^D a_i\left[(k+d)^i-k^i\right]}{a}.$$
    Since $d\mid k$, we have $d^i\mid (k+d)^i-k^i$. 
    Note that for all $q\in \{p_i,q_j \,|\, i=1,\ldots,r,\, \, j=1,\ldots,s\}$ we have $\nu_q(d^2)>\nu_q(a_1d)$, and so $$\nu_q\left(p(k+d)-p(k)\right)=\nu_q(a_1d),$$
    yielding that
    $$ \nu_q(m)=\nu_q(a_1)+\nu_q(d)-\nu_q(a).$$
    Hence, $\nu_q(m)=\nu_q(u)$, so $u\mid m$ and $\gcd(m',2ab)=\gcd(m',u)=1$.
    \item $(a+b)\cdot 0 = 0 \equiv p(0) \pmod{u}$
    \item We have $u\mid p(k)$, since $d\mid k$, $u\mid d^2$ and $u\mid a_1d$. So $(a+b)\cdot 0 = 0 \equiv p(k) \pmod{u}$.
    \item We have $b\mid p(k)-aj$, since $b\mid u$, $u\mid j$ and $u\mid p(k)$ by (5).
\end{enumerate}
\end{proof}

\begin{remark}
This proof shows that the statement is still true if we add a constant term $a_0$ to $p(z)$ that satisfies $d\mid a_0$ and $u\mid a_0$. Here $d\mid a_0$ ensures that (1) remains true and $u\mid a_0$ ensures that (4) and (5) will still hold.

Note that the role of $a$ and $b$ in the proof can be interchanged to get a different value of $u$ that can be used instead in this observation.
\end{remark}

\bigskip

\begin{proof}[Proof of Theorem~\ref{cor-b=1}]
Let us write $n$ as $n=2^{\gamma}\prod\limits_{i=1}^r p_i^{\gamma_i}$  where the $p_i$ are distinct odd primes, $\gamma$ is a nonnegative integer and the $\gamma_i$ are positive integers. Let us express $a$ in the form 
$a=2^{\alpha}a'\prod\limits_{i=1}^r p_i^{\alpha_i}$, where $\alpha$ and the $\alpha_i$ are nonnegative integers and $\gcd(a',2n)=1$. 
Without loss of generality we can assume that $\alpha=\alpha'n-\alpha_0$ for some $\alpha',\alpha_0\in\zz^{\ge 0}$. The options for $\alpha'$ and $\alpha_0$ need not to be unique, however if $\alpha=0$ then we choose $\alpha'=\alpha_0=0$. 
    Now, let
    \[
    d=2^{\alpha'} a'',\text{ and }\,u=2^{\alpha_0}\cdot\prod\limits_{\substack{p\mid \gcd(n,a)\\p\ge 3 \text{ prime}}} p^{\nu_p(n)},
    \]
where $a''=a'\prod\limits_{i=1}^r p_i^{\alpha_i}$. It is easy to see that $u=2^{\alpha_0}\cdot\prod\limits_{i\in J} {p_i}^{\nu_{p_i}(n)}$ is an equivalent definition for $u$, where 
    \[
    J=\set{i\in\set{1,\,2,\dots,r}:p_i\mid a''}=\set{i\in\set{1,\,2,\dots,r}:\alpha_i\neq 0}.
    \]
    
We will define $t$ in such a way that it satisfies the congruence relations
    \[
    \left\lbrace{
    \begin{array}{cl}
        t\equiv 0&\pmod{2^{\alpha'}}   \\
        t\equiv 1&\pmod{a''}
    \end{array}
    }\right.
    \]
and we will define $v$ such that $v\equiv\left(a+1\right)^{-1} t^n\pmod{u}$. This definition is meaningful, because of the following two facts:

\begin{itemize}
    \item First we used that $a+1$ has a multiplicative inverse modulo $u$.
    \begin{proof}
        For every  prime divisor $p$ of $u$, either $p=p_i$ for some $i\in J$, or $p=2$. If $p=p_i$, then $p\mid a''\mid a$, hence $p\nmid a+1$. If $p=2$, then $\alpha_0\neq0$, which implies that $\alpha\neq 0$, so $p\mid a$ and $p\nmid a+1$. We conclude that $\gcd(u,a+1)=1$, from where our claim follows.
    \end{proof}
    
    \item Secondly, since $t$ is a residue class modulo $d=2^{\alpha'}a''$ and not modulo $u$ we need to prove that if  $t_1\equiv t \pmod {d}$ then $\left(a+1\right)^{-1}t_1^n\equiv\left(a+1\right)^{-1}t^n\pmod{u}$. This result will be a consequence of the following claim.
   
\end{itemize}

\bigskip

\textit{Claim.} If $t_1\equiv t\pmod{d}$, then $t_1^n\equiv t^n\pmod{u}$, where $d,\,u,\,t$ are as defined before.

\smallskip

\textit{Proof.} Indeed if $t_1\equiv t \pmod {d}$, then $t_1=t+t_2d$ for some integer $t_2$. From the definition of $t$ we get that $t=2^{\alpha'}t'$ for some integer $t'$ , hence
    \[
    t_1^{n}-t^n=\left(2^{\alpha'}t' +t_2 2^{\alpha'}a''\right)^n-\left(2^{\alpha'}t'\right)^n
    =2^{\alpha' n}\left[
    \left(t' +t_2 a''\right)^n-\left(t'\right)^n
    \right],
    \]
so $\nu_2(t_1^n-t^n)\geq \alpha' n=\alpha+\alpha_0\geq\alpha_0=\nu_2(u)$.
If we take another prime $p\geq 3$ such that $p\mid u$, then we know that $p=p_i$ for some $i\in J$. As $p>2$, we have $p\mid a''$, which implies
    \[
    \gcd(p,t)=\gcd(p,t_1)=1,
    \]
    because $t$ and $t_1$ are invertible modulo $a''$.
     We also have $p=p_i\mid d t_2=t_1-t$, so we can apply the so-called lifting-the-exponent lemma which gives us
\begin{align*}
    \nu_p\left(t^n_1-t^n\right)
    &=\nu_p(n)+\nu_p(d t_2)=\gamma_i+\nu_p(2^{\alpha'})+\nu_p(a'')+\nu_p(t_2) \\
    &= \gamma_i+0+\alpha_i+\nu_p(t_2)\geq \gamma_i=\nu_p(u).
\end{align*}
Since $\nu_p(u)\leq \nu_p(t_1^n-t^n)$ for every prime divisor $p$ of $u$, we have $t_1^n\equiv t^n\pmod{u}$. $\qed$

\bigskip

Now, we are ready to prove that the conditions of Theorem~\ref{mainresult} are satisfied.
Indeed:
\begin{enumerate}
   
    \item The relation $(a+1)t\equiv t^n$ is satisfied modulo $2^{\alpha'}$ and modulo $a''$, furthermore, $\gcd(2^{\alpha'},a'')=1$, hence it is also satisfied modulo $2^{\alpha'}a''=d$;
   
    \item Let $k\equiv t \pmod{d}$. We know that $k+d\equiv k\equiv t\pmod{d}$, so using the reasoning in the proof of the claim above, for every $i\in J$, $p=p_i$, we have
    \[
    \nu_p((k+d)^n-k^n)=\gamma_i+\alpha_i+\nu_p(1)=\gamma_i+\alpha_i=\nu_p(au)\geq \alpha_i=\nu_p(a),
    \]
    and for $p=2$ we have
    \[
    \nu_p((k+d)^n-k^n)=\alpha+\alpha_0=\nu_p(au)\geq \alpha=\nu_p(a).
    \]
    Now, we will prove that $a'\mid (k+d)^n-k^n$, so then together with our previous statements we have that $a\mid (k+d)^n-k^n$.
    Indeed, because of $k\equiv t \equiv 1 \pmod{a''}$, we have that $k=1+k_1 a'$ for some integer $k_1$, and also $d=d_1a'$, where $d_1=2^{\alpha'}\prod\limits_{i=1}^r p_i^{\alpha_i}$. Using the lifting-the-exponent lemma, for every prime $p\mid a'$, we have that
\begin{align*}
 \nu_p((k+d)^n-k^n)
 &=\nu_p\left( (1+\left(k_1+d_1)a'\right)^n-(1+k_1a')^n \right)\\
 &=\nu_p(n)+\nu_p(a' d_1)=\nu_p(a'),
\end{align*}
where the last equation follows from  $\gcd(a',n)=\gcd(a',d_1)=1$.

    \item From the proof of $(2)$ it is easy to see that $au\mid (k+d)^n-k^n$, hence $u\mid m$.
    Once again using the results in $(2)$, and the fact that $\nu_p(m')=\nu_p((k+d)^n-k^n)-\nu_p(u)-\nu_p(a)$ for every prime $p$, we can write  that 
    \[
    \nu_p(m')=
    \left\lbrace{
    \begin{array}{ll}
        \alpha+\alpha_0-\alpha_0-\alpha, & \text{if } p=2;  \\
        \gamma_i+\alpha_i-\gamma_i-\alpha_i, & \text{if }p=p_i,\,i\in J;\\
        \nu_p(a')-\nu_p(a'), & \text{if }p\mid a';\\
        \nu_p(m'), &\text{ otherwise},
    \end{array}}\right. 
    \]
    from where the rest of the conditions in (3) follow.
   
    \item From $(1)$ an application of our Claim yields that $t^n\equiv\left(\left( 1+a \right)^{-1}t^n \right)^n \pmod{u}$, which is equivalent with
    $(1+a)v\equiv v^n\pmod{u}$.
   
    \item Using the fact that $(1+a)v\equiv(1+a)(1+a)^{-1} t^n\pmod{u}$, we can once again apply our Claim to get that $k^n\equiv (1+a)v\pmod{u}$.
   
    \item From $b=1$ it follows that $b\mid p(k)-aj$, for every $k,j\in\mathbb{Z}$.
\end{enumerate}
\end{proof}

\bigskip

\bigskip

\begin{proof}[Proof of Corollary~\ref{cor-czp}]
Let $d=aBp^n$ and $t$ be as in the statement of the theorem.

Let us choose $u$ as follows. 
$$u=\begin{cases}
2^{\nu_2(c)}p^{\nu_p(c)+pn}B, & \text{if }p\neq2\\
p^{\nu_p(c)+pn}B, & \text{if }p=2
\end{cases}$$
Let us pick $v\in \zz_u$ such that the congruences
\[
    \left\lbrace{
    \begin{array}{ll}
        v \equiv 0 &\pmod{p^{\nu_p(c)+pn}}\\
        v \equiv 0 &\pmod{2^{\nu_2(c)}}\\
    v \equiv t &\pmod{B}
    \end{array}
    }\right.\]
hold. 
Note that such a $v$ exists in all cases by the assumptions of the theorem (notice that if $2\mid B$, then $\nu_2(c)=0$). These choices satisfy the conditions of Theorem~\ref{mainresult}, as if we assume $k\equiv t \pmod{d}$ and $j\equiv v \pmod{u}$:

\begin{enumerate}
    \item $(a+b)\cdot t \equiv ct^p \pmod{d}$ holds by the assumptions on $t$.
    \item Since $a\mid d$, we have $a\mid c(k+d)^p-ck^p$.
    \item First note that since $a$, $B$ and $c$ are pairwise coprime and $\gcd(a,b)=1$, $t$ is coprime to both $a$ and $B$ and hence so is $k$.
    \noindent
    We have $$m=c\cdot\frac{(k+d)^p-k^p}{a}=c\cdot\sum_{i=1}^p \binom{p}{i}\frac{d^i}{a}k^{p-i}.$$ As $a\mid d$, every term in the sum is divisible by $a$, except for the first term, $\frac{pdk^{p-1}}{a}$, which is coprime to $a$ using the assumption $\gcd(a,p)=1$. Hence $\gcd(a,m)=1$. Similarly, $B\mid m$, but $\gcd\left(\frac{m}{B},B\right)=1$.
    
    As $p^n\mid k$, write $k=p^n\ell$ for some $\ell\in\zz^+$. Then $$(k+d)^p-k^p=p^{pn}\left((\ell+aB)^p-\ell^p\right).$$
     As $(\ell+aB)^p-\ell^p\equiv \ell+aB-\ell=aB \not\equiv 0 \pmod{p}$, we have $\nu_p(m)=\nu_p(c)+pn$.
    
    Hence, if $p=2$, then $m'=\frac{m}{u}\in \zz$ and $\gcd(m',a)=\gcd(m',B)=\gcd(m',2)=1$.
    
    If $p\neq 2$, then we also have $m''=\frac{m}{p^{\nu_p(c)+pn}B}\in\zz$ that satisfies $\gcd(m'',a)=\gcd(m'',B)=\gcd(m'',p)=1$. If $2\mid a$ or $2\mid B$, then $\nu_2(c)=0$, $m''=m'$ and $\gcd(m',2)=1$. Otherwise, $m'=\frac{m''}{2^{\nu_2(c)}}=\frac{m}{u}\in \zz$ and $\gcd(m',2)=1$, since $2\nmid (k+d)^p-k^p$, as $d$ is odd.
    
    Hence, in both cases $u\mid m$, $2\nmid m'$ and $(a,m')=(b,m')=(u,m')=1$.
    
    \item $(a+b)v \equiv 0 \equiv cv^p \pmod{2^{\nu_2(c)}}$, $(a+b)v\equiv 0 \equiv cv^p \pmod{p^{\nu_p(c)+pn}}$ and $(a+b)v \equiv at \equiv ct^p \equiv cv^p \pmod{B}$, by the choice of $t$. So the equivalence holds modulo $u$, as well.
    \item $(a+b)v\equiv 0 \equiv ck^p \pmod{2^{\nu_2(c)}}$ and $(a+b)v\equiv 0 \equiv ck^p \pmod{p^{\nu_p(c)+pn}}$, since $p^n\mid k$ and $(a+b)v \equiv at \equiv ct^p \equiv ck^p \pmod{B}$. So the equivalence holds modulo $u$, as well.
    \item We have $B\mid ck^p-aj$ from (5), as $aj\equiv (a+b)v \pmod{B}$. Also,  $p^{pn-\varepsilon}\mid ck^p-aj$, since $p^{pn}$ divides both $k^p$ and $j$. Therefore, $b\mid ck^p-aj$.
\end{enumerate}
\end{proof}

\bigskip 

\begin{proof}[Proof of Corollary~\ref{rem-cz2}]
Without loss of generality, we may assume that $a$ is odd. Now, Corollary~\ref{cor-czp} may be applied for the equation $ax+by=cz^2$, since the condition about the solvability of the system of congruences automatically holds when $p=2$. Hence, $ax+by=cz^2$ is 2-Ramsey.
\end{proof}

\bigskip

\noindent
The following lemma can be used to extend the 2-Ramseyness of the equation $ax+by=cz^2$ (where $a,b,c$ are pairwise coprime) to $d(ax+by)=cz^2$.

\begin{lemma}
\label{general}
Assume that the equation $ax+by=c^{-2}p(cz)=\sum\limits_{i=0}^D a_ic^{i-2}z^i$ is 2-Ramsey, where $a,b,c,a_D\in \zz^+$, $p\in\zz[z]$, $c\mid a_1$ and $c^2\mid a_0$. Then the equation $acx+bcy=p(z)=\sum\limits_{i=0}^D a_iz^i$ is 2-Ramsey.
\end{lemma}

\begin{proof}
Let $\phi$ be a 2-coloring of $\zz^{+}$. Let us consider only those solutions of the equation $acx+bcy=p(z)$ when $x,y,z$ are all divisible by $c$. Let $x=cx'$, $y=cy'$ and $z=cz'$. Substituting into the equation, the triple $x,y,z$ is a solution if and only if $$ax'+by'=\sum_{i=0}^D a_ic^{i-2}\left(z'\right)^i=c^{-2}p(cz').$$ Let $\psi(j):=\phi(cj)$ be another 2-colouring of $\zz^{+}$. Then by the assumption, there are infinitely many solutions $x',y',z'$ such that $ax'+by'=c^{-2}p(cz')$ and $\psi(x')=\psi(y')=\psi(z')$. But then for such solutions we have $acx+bcy=p(z)$ and by the definition of $\psi$, we have $\phi(x)=\phi(y)=\phi(z)$. Hence, the equation $acx+bcy=p(z)$ is 2-Ramsey.
\end{proof}

\bigskip

\begin{proof}[Proof of Corollary~\ref{cor-dcz2}]
The statement follows directly from Corollary~\ref{rem-cz2} and Lemma~\ref{general}.
\end{proof}

\section{Concluding remarks}\label{sec-remark}

In this paper we prove that the equation $ax+by=p(z)$ is 2-Ramsey if the coefficients $a,b$ and the polynomial $p$ satisfy certain conditions. On the other hand, in all the cases when we can prove that 2-Ramseyness fails to hold there is a 2-colouring avoiding monochromatic solutions of the following type. Let $m$ be a positive integer and $\mathbb{Z}_m=A\cup B$ be a partitioning of the residues in such a way that $(aA+bA)\cap p(A)=(aB+bB)\cap p(B)=\emptyset$. If we colour the integers having a modulo $m$ residue lying in $A$ to the first colour and lying in $B$ to the second colour, then there will be no monochromatic solution. Note that is a periodic colouring. It would be interesting to decide whether a periodic colouring avoiding monochromatic solutions always exists when $ax+by=p(z)$ is not 2-Ramsey. Note that in the case $a=b=1$ when $x+y=p(z)$ is not 2-Ramsey, then in fact {\it every} 2-colouring avoiding monochromatic solutions is periodic according to \cite{LPS}.

\section{Acknowledgements}
We would like to thank the anonymous referees for their useful comments and suggestions. The research was supported by the Lend\"ulet program of the Hungarian Academy of Sciences (MTA). PPP was also supported by the National Research, Development and Innovation Office NKFIH (Grant Nr. K124171 and K129335).

\end{document}